# Nonlinear Control Synchronization Method for Fractional-order Time Derivatives Chaotic Systems

**Vivek Mishra[a], S. K. Agrawal[b]**

[a]Alliance School of Applied Mathematics, Alliance University, Bengaluru-562106, India.
[b]Department of Applied science, Bharati Vidyapeeth College of Engineering, New Delhi-110063, India

**Abstract**: "Synchronization of two dynamical systems" is the term used to describe the phenomenon when two or more systems gradually change their states or behaviors to become similar or identical. This can happen in a lot of fields, such as physics, engineering, biology, and economics. Synchronization finds applications in neurology and communication systems. It is present in both man-made and organic systems. The nonlinear control synchronization technique for fractional-order time derivative systems is described in this article, where the Adams Basford Moulton method is used for solving the fractional-order system. The reliability and ease of applicability for two chaotic systems are demonstrated by the numerical simulation. Furthermore, in this article, both systems were kept in a chaotic condition while being synchronized with each other. The effects of synchronizing time and rearranging the derivatives are the most significant sections of this article.

**Keywords**: Fractional time derivative; Synchronization; Financial system; Volta's chaotic system; Chaotic system.

## 1. Introduction

Currently, nonlinear phenomena in different fields of science and engineering have attracted a lot of attention from scientists and engineers. In the last few years, they have worked hard to develop models that work well by using nonlinear dynamics. Fractional-order derivatives and integrals in nonlinear models have added new dimensions to existing problems. In the field of nonlinear sciences, theories and results in chaotic dynamics, chaos control, synchronization, and the boundedness of chaotic attractors have opened a new area of study. Most of the nonlinear systems demonstrate chaos, which sensitively depends on initial conditions. When many chaotic systems modify some aspect of their motion as a shared behavior, this is known as chaos synchronization. It is because of forcing (periodical or noisy) or coupling.

The characteristics of fractional order derivatives for which they have gained much popularity are the nonlocality and memory effect. The integer order derivative is the local operator, whereas the fractional-order derivative uses all nonlocal values of the function. The fractional-order chaotic dynamical system has been investigated and studied in the fields of physical and mathematical communities in the last few decades. Mathematically, synchronization is achieved when

$\lim_{t\to\infty}\|x_1(t)-x_2(t)\|=0$, where $x_1(t)$ and $x_2(t)$ is respectively the state vectors of the drive and response systems.

The pioneering work of Pecora and Carroll [20], introduced the feasibility of synchronizing chaotic systems through simple coupling with different initial conditions. It has attracted a huge amount of interest among scientists from different fields because of its significant applications in the physical system [14], chemical system [10], ecological system [1], brain activity models, system identification, secure communication, and pattern recognition phenomena ([6,13]), etc. In recent years, many synchronization schemes have been proposed, such as linear and nonlinear feedback synchronization [3,12,30], time delay feedback approach [9,19], active control [11,26], adaptive control [4,23], sliding mode control [5,28], projective synchronization [27,29], phase synchronization [18], etc. Currently, nonlinear control has been applied to synchronize between identical and nonidentical systems. This method is mostly used with the stability criterion of the nonlinear fractional-order differential equation, where the considered Jacobian matrix A of the vector field satisfies $|\arg(spec(A))| > q\pi/2$ [17] or those eigenvalues of matrix A satisfying the same condition have geometric multiplicity one.

In this article, the authors have investigated the synchronization between fractional-order financial and Volta's chaotic systems by the nonlinear control method [7, 8]. Numerical simulations are carried out for different fractional-order derivatives, which are displayed graphically to express the efficiency of the proposed approach.

## 2. System descriptions

### 2.1 The fractional-order Financial System [2]

$$\frac{d^{q_1} x_1}{dt^{q_1}} = z_1 + (y_1 - \alpha)x_1$$

$$\frac{d^{q_2} y_1}{dt^{q_2}} = 1 - \beta y_1 - x_1^2 \qquad (1)$$

$$\frac{d^{q_3} z_1}{dt^{q_3}} = -x_1 - \gamma z_1$$

Here, $\alpha$, $\beta$ and $\gamma$ are the saving amount, the cost per investment, and the elasticity of demand in the commercial market, respectively. The state variables are $x_1(t)$ is the interest rate, $y_1(t)$ is the investment demand, and $z_1(t)$ is the price index. During synchronization, the values of parameters are initially taken as [2, -1, 1] then the system (1) shows chaotic behavior with the derivative's order $q_i \geq 0.8436$. The chaotic attractors of the system (1) are described in Fig.1 for the order of derivatives $q_i = 0.99$.

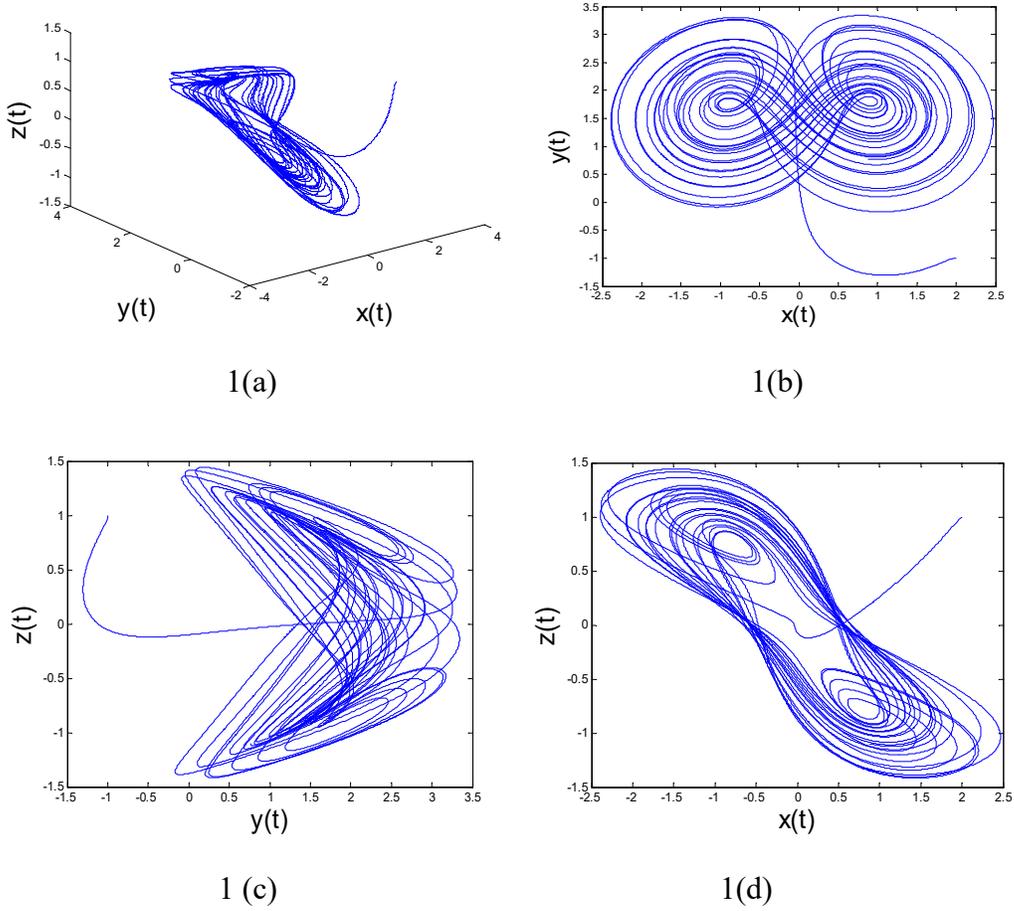

1(a)

1(b)

1 (c)

1(d)

**Fig.1.** Phase portraits of the financial system.

## 2.2 The fractional-order Volta's System [21]

The fractional-order Volta's System is given by

$$\frac{d^{q_1} x_2}{dt^{q_1}} = -x_2 - ay_2 - z_2 y_2$$
$$\frac{d^{q_2} y_2}{dt^{q_2}} = -y_2 - bx_2 - x_2 z_2$$
(2)

$$\frac{d^{q_3} z_2}{dt^{q_3}} = cz_2 + x_2 y_2 + 1$$

Here, a, b, and c are variable parameters. For the synchronization, the parameters, values are taken as a=19, b=11, and c=.73, with the initial condition [8, 2, 3]. In this case, the system (2) can express chaotic behavior with commensurate order of derivatives as $q_i = 0.98 (i = 1,2,3)$. the chaotic attractors of the system (2) are described in Fig.2 for the order of the derivative $q_i = 0.99$.

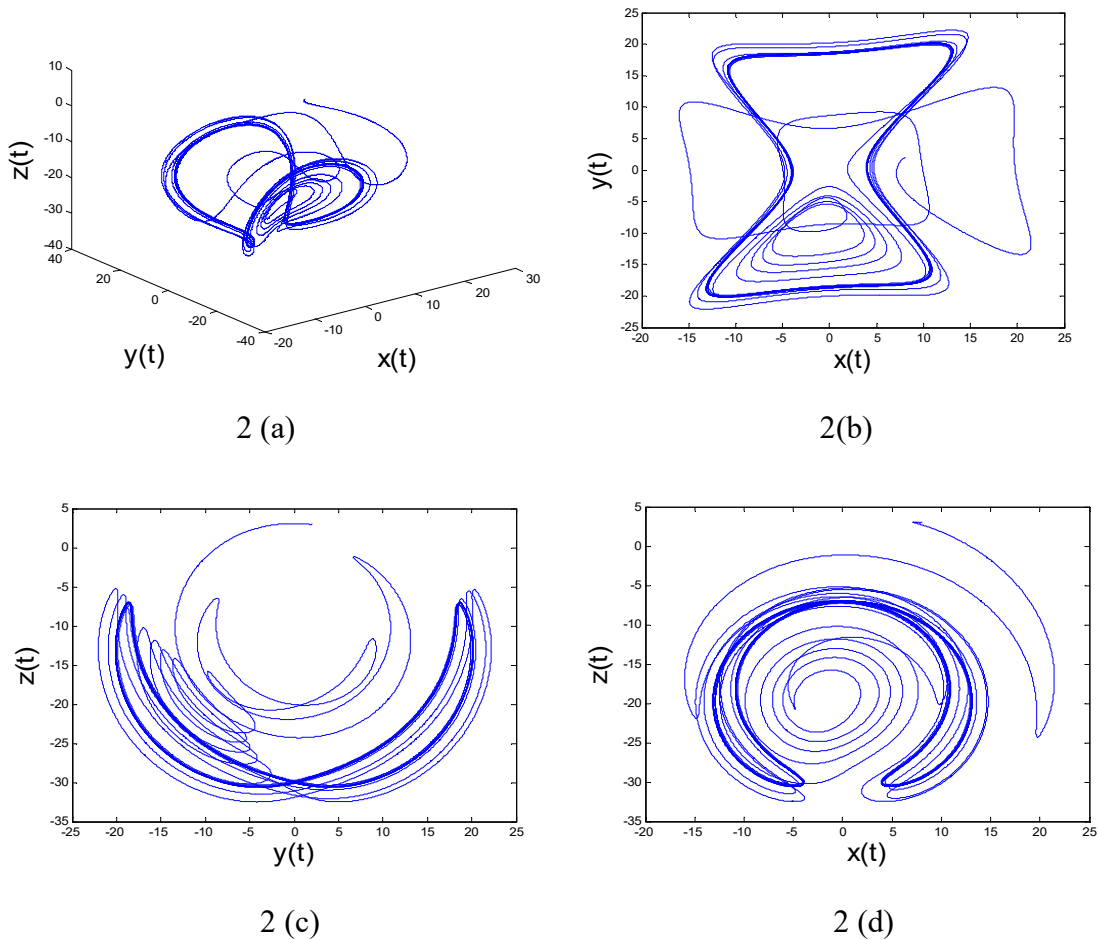

2 (a)

2(b)

2 (c)

2 (d)

**Fig.2.** Phase portraits of Volta's system.

### 3.1. Synchronization of fractional-order Financial and fractional-order Volta's systems via nonlinear control method.

The synchronization behavior between two different fractional orders Financial and Volta's system is achieved in this section. Let the Financial system drives Volta's system; therefore, the financial system is taken as the master system and Volta's as a slave system. The synchronization behavior between two different fractional orders in the financial and Volta systems is achieved in this section. Let the financial system drive Volta's system; therefore, the financial system is taken as the master system and Volta's as a slave system.

The master system is given by (1)

$$\begin{cases} \frac{d^{q_1} x_1}{dt^{q_1}} = z_1 + (y_1 - \alpha)x_1 \\ \frac{d^{q_2} y_1}{dt^{q_2}} = 1 - \beta y_1 - x_1^2 \\ \frac{d^{q_3} z_1}{dt^{q_3}} = -x_1 - \gamma z_1 \end{cases} \quad (3)$$

The slave system is described by (2)

$$\begin{cases} \frac{d^{q_1} x_2}{dt^{q_1}} = -x_2 - ay_2 - z_2 y_2 + u_1(t) \\ \frac{d^{q_2} y_2}{dt^{q_2}} = -y_2 - bx_2 - x_2 z_2 + u_2(t) \\ \frac{d^{q_3} z_2}{dt^{q_3}} = cz_2 + x_2 y_2 + 1 + u_3(t), \end{cases}$$
(4)

Here $u_1(t), u_2(t)$ and $u_3(t)$ are nonlinear control functions. We aim to analyze the synchronization of systems (1) and (2). We define the error states as $e_1 = x_2 - x_1, e_2 = y_2 - y_1, e_3 = z_2 - z_1$. Hence, the error dynamics are obtained by equation (2) and equation (4) as

$$\begin{cases} \frac{d^{q_1} e_1}{dt^{q_1}} = -e_1 - ae_2 + e_3 + (\alpha - 1)x_1 - (x_1 + a)y_1 - (1 + y_2) + u_1(t) \\ \frac{d^{q_2} e_2}{dt^{q_2}} = -be_1 - e_2 + (\beta - 1)y_1 - (b - x_1)x_1 - x_2 z_2 - 1 + u_2(t) \\ \frac{d^{q_3} e_3}{dt^{q_3}} = -e_1 + ce_3 + (y_2 + 1)x_2 + (c + \gamma)z_1 + 1 + u_3(t) \end{cases}$$
(5)

The nonlinear control inputs $u_1(t), u_2(t), u_3(t)$ are defined as

$$\begin{cases} u_1(t) = -(\alpha - 1)x_1 + (x_1 + a)y_1 + (1 + y_2) + v_1(t) \\ u_2(t) = -(\beta - 1)y_1 + (b - x_1)x_1 + x_2 z_2 + 1 + v_2(t) \\ u_3(t) = -(y_2 + 1)x_2 - (c + \gamma)z_1 - 1 + v_3(t) \end{cases}$$
(6)

which leads to

$$\begin{cases} \dfrac{d^{q_1}e_1}{dt^{q_1}} = -e_1 - ae_2 + e_3 + v_1(t) \\ \dfrac{d^{q_2}e_2}{dt^{q_2}} = -be_1 - e_2 + v_2(t) \\ \dfrac{d^{q_3}e_3}{dt^{q_3}} = -e_1 + ce_3 + v_3(t) \end{cases} \qquad (7)$$

The synchronization error of the system in (7) is a linear system with nonlinear control inputs, $V_1(t), V_2(t)$ and $V_3(t)$. Next to frame a feedback control that stabilizes the system such that $e_1, e_2, e_3$ converge to zero as time t goes to infinity, which means that Financial and Volta's systems are synchronized with feedback control. There are many possible options for the control inputs $V_1(t), V_2(t)$ and $V_3(t)$. We have framed

$$\begin{bmatrix} v_1(t) \\ v_2(t) \\ v_3(t) \end{bmatrix} = A \begin{bmatrix} e_1 \\ e_2 \\ e_3 \end{bmatrix} \qquad (8)$$

where A is a matrix of constants of order 3x3. For stability of the closed loop system, matrix A should be selected in such a way that the feedback system has eigenvalues $\lambda_i$ that satisfies the control $|arg(\lambda_i)| > 0.5\pi q, i = 1,2,3$. Matrix A is not unique. It can be chosen as follows:

$$A = \begin{pmatrix} 0 & a & -1 \\ b & 0 & 0 \\ 1 & 0 & -1-c \end{pmatrix} \qquad (9)$$

Then the error system is changed to

$$D_t^{q_1} e_1 = -e_1, \; D_t^{q_2} e_2 = -e_2, \; D_t^{q_3} e_3 = -e_3. \qquad (10)$$

Here all eigenvalues $\lambda_i$ of the matrix are -1, which satisfies the condition $|arg(\lambda_i)| > q\pi/2$, for $0 < q \leq 1$. Hence, the linear system (10) is stable, and we achieved the required synchronization.

### 3.2 Numerical simulation and results

For numerical simulations, the parameters of the Financial and Volta's systems are considered as $\alpha = 1, \beta = 0.1, \gamma = 1$, and $a = 19, b = 11, c = 0.73$ respectively, and the step size for time variation is considered as 0.0005. The initial values of the drive system and the response system are $(2, -1, 1)$ and, (8,2,3) respectively. Hence, the initial errors are $(6,3,2)$. Now choose $\lambda_1 = -1, \lambda_2 = -1, \lambda_3 = -1$, the control functions that can be evaluated and the synchronization is achieved between the drive and response system. State trajectories of the synchronization between the systems are depicted in Fig.3 and Fig.4 for the order of the derivative, $q_i = 0.99$ and $q_i = 0.98, i = 1,2,3$ respectively.

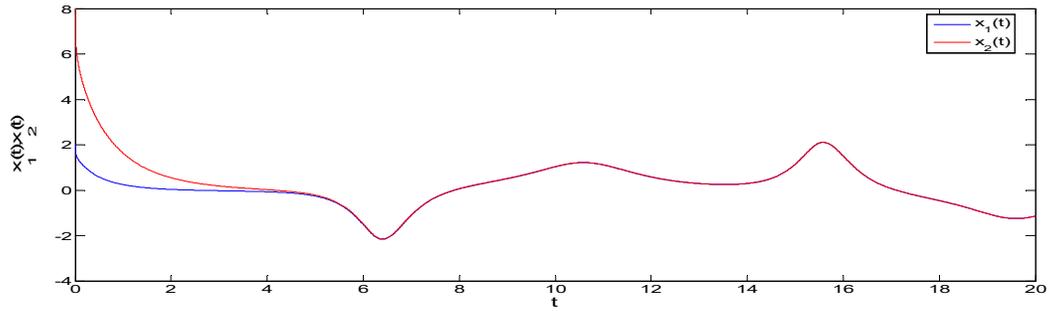

3(a)

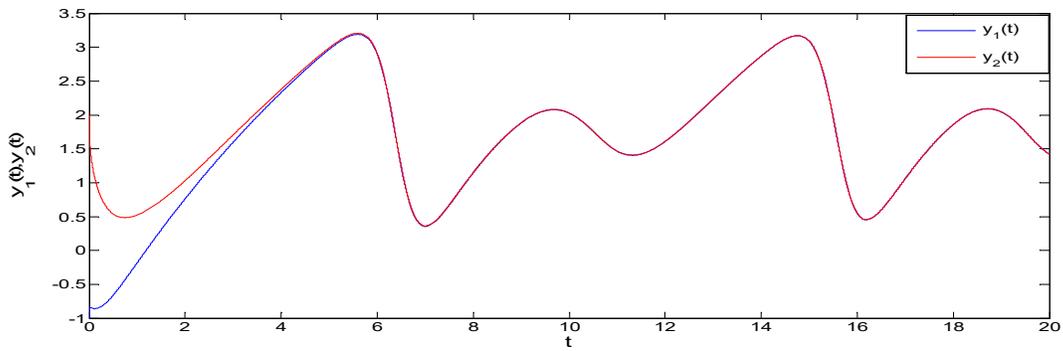

3(b)

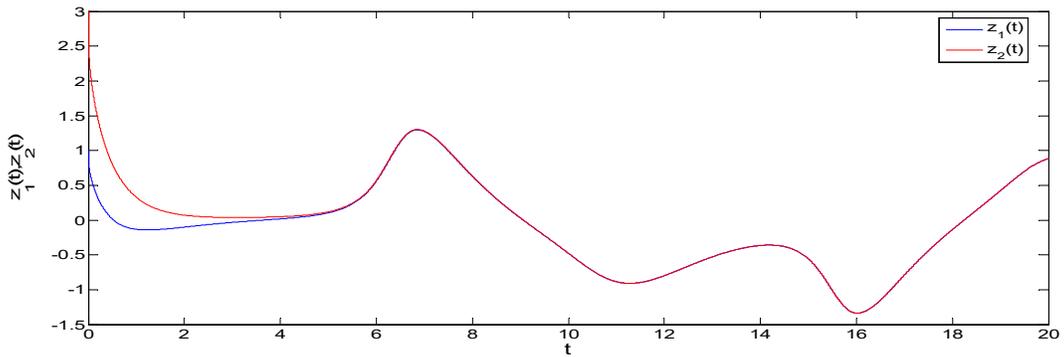

3(c)

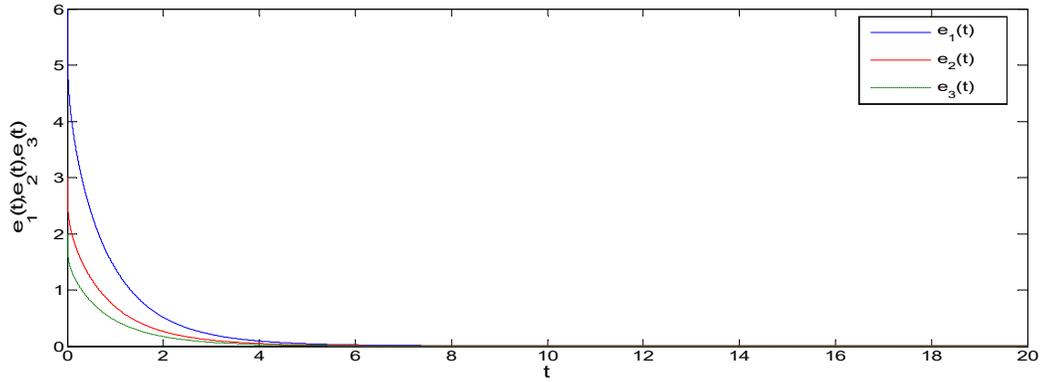

3(d)

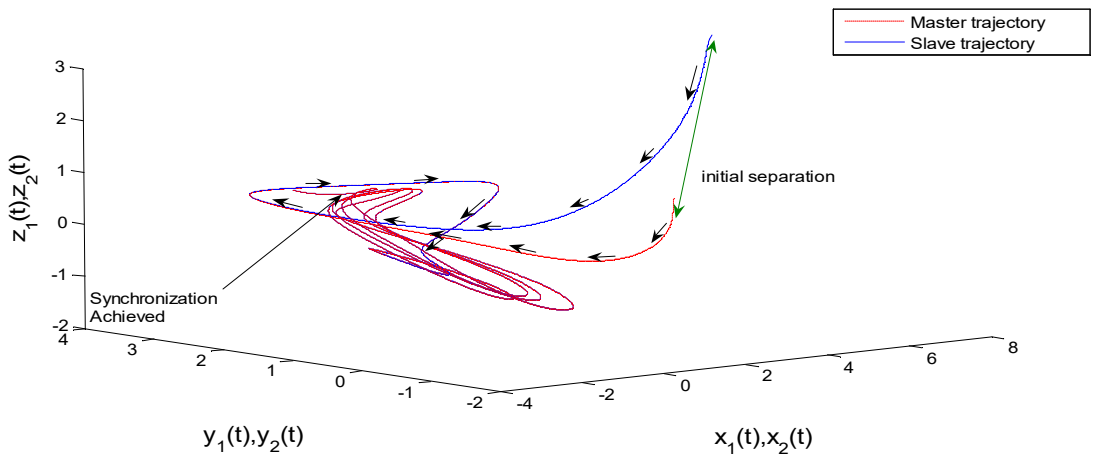

3(e)

**Fig.3.** State trajectories of the drive system (3) and response system (4); Fig.3 (a): between $x_1$ and $x_2$, Fig. 3(b): between $y_1$ and $y_2$, Fig. 3(c): between $z_1$ and $z_2$, Fig. 3(d): The evolution of the error functions e1 (t), e2 (t) and e3 (t), (t), Fig. 3(e): Phase space representation of the synchronization of the master and slave chaotic systems.

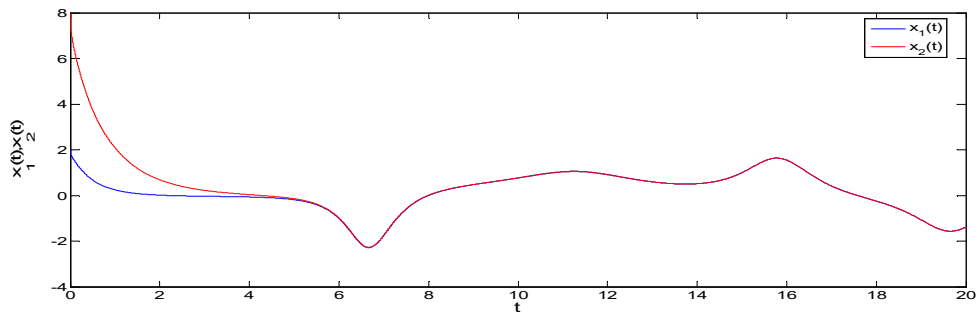

4(a)

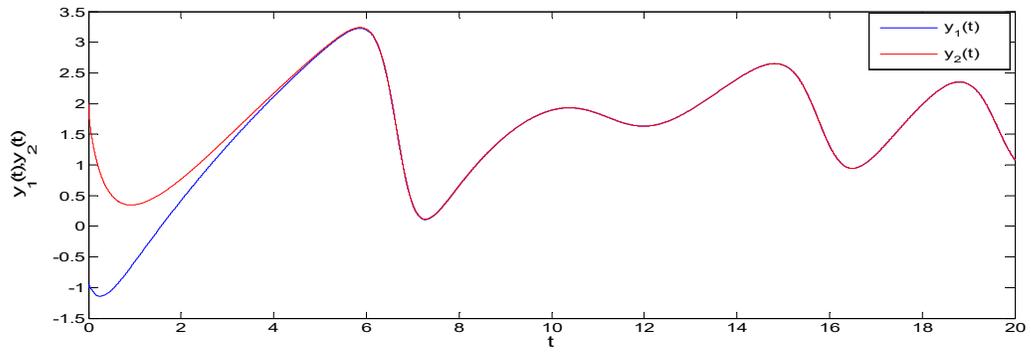

4(b)

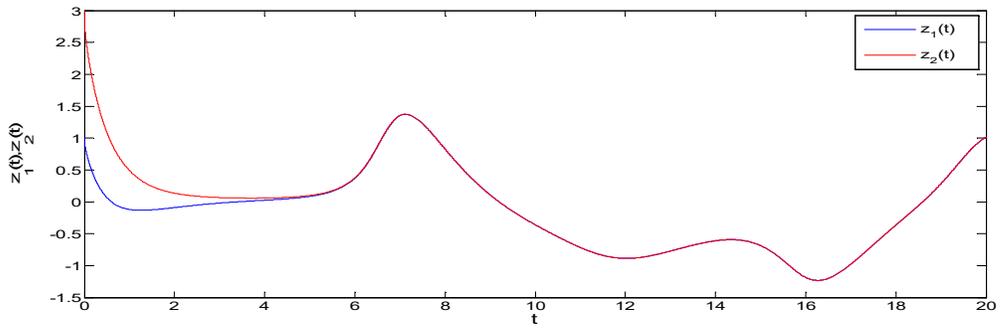

4(c)

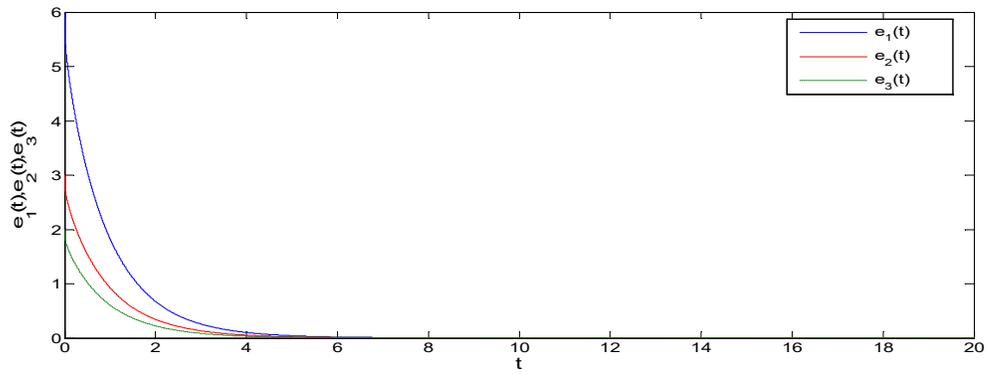

4(d)

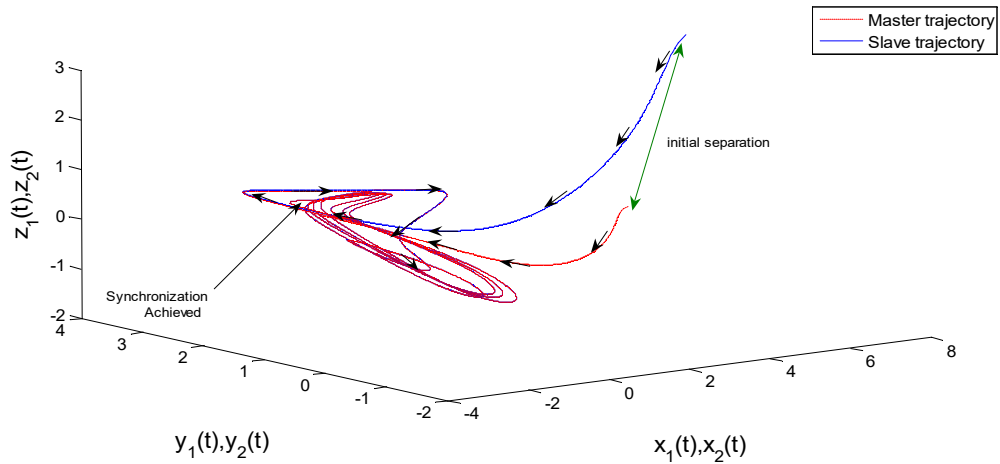

4(e)

**Fig.4.** Trajectories of the states of drive system (3) and response system (4); Fig.4 (a): between $x_1$ and $x_2$, Fig. 4(b): between $y_1$ and $y_2$, Fig. 4(c): between $z_1$ and $z_2$, Fig. 4(d): The evolution of the error functions $e_1(t)$, $e_2(t)$ and $e_3(t)$, Fig. 4(e): Phase space representation of the synchronization of the master-slave chaotic systems.

Fig.3 and Fig.4 represent that the systems are synchronized after a time for both the fractional-order time derivative. It is observed from the above figures that the time of synchronization decreases with a decrease in the order of derivatives. Thus, the synchronization time decreases as the systems approach their stability. Graphs are plotted using MATLAB 2007.

## 4. Conclusion

This study aims to present the synchronization between two different chaotic systems of fractional order. Considering the stability analysis, the synchronization of chaotic systems of fractional order through linear controller input parameters for each system has been achieved. The design of the controller that is applied to the response system affects the system dynamics to achieve synchronization. In the numerical simulation result, the design of the controller is such that the components of the error system reach zero as the time goes to infinity, which demonstrates that the method is easy to implement for the synchronization of two chaotic systems of fractional order. Another important objective of the study is to show the variations in the time of synchronization as the systems approach stability with a decrease in fractional-order derivatives.